\documentstyle [12pt]{article}

\def\qmod#1#2{{\hbox{}^{\displaystyle{#1}}}\!\big/\!\hbox{}_{
\displaystyle{#2}}}

\newfam\msbfam
\font\twelmsb=msbm10 at 12pt
\font\tenmsb=msbm10
\font\sevenmsb=msbm10 at 7pt
\font\fivemsb=msbm10 at 5pt

\textfont\msbfam=\twelmsb
\scriptfont\msbfam=\sevenmsb
\scriptscriptfont\msbfam=\fivemsb

\def\Bbb{\fam\msbfam\tenmsb}

\def\C{{\Bbb C}}

\def\R{{\Bbb R}}
\def\Z{{\Bbb Z}}

\def\AA{{\cal A}}
\def\BB{{\cal B}}
\def\CC{{\cal C}}
\def\DD{{\cal D}}
\def\EE{{\cal E}}
\def\FF{{\cal F}}

\def\MM{{\cal M}}

\def\OO{{\cal O}}

\def\SS{{\cal S}}

\def\VV{{\cal V}}
\def\WW{{\cal W}}

\def\qed {\hfill\vrule height6pt width6pt depth0pt \bigskip}

\def\map{\longrightarrow}
\def\textmap#1{\mathop{\vbox{\ialign{
                                ##\crcr
    ${\scriptstyle\hfil\;\;#1\;\;\hfil}$\crcr
    \noalign{\kern-1pt\nointerlineskip}
    \rightarrowfill\crcr}}\;}}
\newcommand{\Cal}{\cal}
\def\textlmap#1{\mathop{\vbox{\ialign{
                                ##\crcr
    ${\scriptstyle\hfil\;\;#1\;\;\hfil}$\crcr
    \noalign{\kern-1pt\nointerlineskip}
    \leftarrowfill\crcr}}\;}}

\newfam\meuffam
\font\tenmeuf=eufm10
\font\sevenmeuf=eufm7
\font\fivemeuf=eufm5

\textfont\meuffam=\tenmeuf
\scriptfont\meuffam=\sevenmeuf
\scriptscriptfont\meuffam=\fivemeuf

\def\germ{\fam\meuffam\tenmeuf}

\def\h{{\germ h}}

\def\p{{\germ p}}

\def\s{{\germ s}}

\begin{document}
\def\Pr{{\rm Pr}}
\def\tr{{\rm Tr}}
\def\End{{\rm End}}
\def\Spin{{\rm Spin}}
\def\U{{\rm U}}
\def\SU{{\rm SU}}
\def\SO{{\rm SO}}
\def\PU{{\rm PU}}
\def\spin{{\rm spin}}
\def\u{{\rm u}}
\def\su{{\rm su}}
\def\so{{\rm so}}
\def\pu{{\rm pu}}
\def\Pic{{\rm Pic}}
\def\NS{{\rm NS}}
\def\deg{{\rm deg}}
\def\Hom{{\rm Hom}}
\def\Herm{{\rm Herm}}
\def\Vol{{\rm Vol}}
\def\pf{{\bf Proof: }}
\def\id{{\rm id}}
\def\i{{\germ i}}
\def\im{{\rm im}}
\def\rk{{\rm rk}}
\def\ad{{\rm ad}}
\def\h{{\bf H}}
\def\coker{{\rm coker}}
\def\dv{\bar\partial}
\def\dva{\bar\partial_A}
\def\da{\partial_A}
\def\p{\partial\bar\partial}
\def\pa{\partial_A\bar\partial_A}
\def\Dr{\hskip 4pt{\not}{D}}
\newtheorem{sz}{Satz}[section]
\newtheorem{th}[sz]{Theorem}
\newtheorem{pr}[sz]{Proposition}
\newtheorem{re}[sz]{Remark}
\newtheorem{co}[sz]{Corollary}
\newtheorem{dt}[sz]{Definition}
\newtheorem{lm}[sz]{Lemma}
\title{The Coupled Seiberg-Witten Equations, Vortices, and Moduli Spaces of
Stable Pairs}
\author{Christian Okonek$^*$\and Andrei Teleman\thanks{Partially supported
by: AGE-Algebraic Geometry in Europe, contract No ERBCHRXCT940557 (BBW
93.0187),
and by  SNF, nr. 21-36111.92} }
\date{January, 13-th 1995 }
\maketitle
\setcounter{section}{-1}
\section{Introduction}

Recently, Seiberg and Witten [W] introduced new invariants of 4-manifolds,
which are defined by counting solutions of a certain non-linear differential
equation.

The new invariants are expected to be equivalent to Donaldson's
polyno\-mial-invariants---at least for manifolds of simple type [KM
1]---and they
have already found important applications, like e.g. in the proof of the Thom
conjecture by Kronheimer and Mrowka [KM 2].

Nevertheless, the equations themselves remain somewhat mysterious, especially
from a mathematical point of view.

The present paper contains our attempt to understand and to generalize the
Seiberg-Witten equations by coupling them to connections in unitary vector
bundles,
and to relate their solutions to more familiar objects, namely stable pairs.

Fix a $\Spin^c$-structure on a Riemannian 4-manifold $X$, and denote by
$\Sigma^{\pm}$ the associated spinor bundles. The equations which we will study
are:
$$
\left\{\begin{array}{lcc}
\Dr_{A,b}\Psi &=&0\\
\Gamma(F_{A,b}^+)&=&(\Psi\bar\Psi)_0\end{array}\right.$$
This is a system of equations for a pair $(A,\Psi)$ consisting of a unitary
connection in a unitary bundle $E$ over  $X$, and a
positive spinor $\Psi\in A^0(\Sigma^+\otimes E)$. The symbol $b$ denotes a
connection in the determinant line bundle of the spinor bundles
$\Sigma^{\pm}$  and
\hbox{$\Dr_{A,b}:\Sigma^+\otimes E\map\Sigma^-\otimes E$} is the Dirac operator
obtained by coupling the connection in $\Sigma^+$ defined by $b$ (and by the
Levi-Civita connection in the tangent bundle)  with the variable connection
$A$ in $E$.

These equations specialize to the original Seiberg-Witten equations
if
$E$ is a line bundle. We show that the coupled equations can be interpreted
as a
differential version of the generalized vortex equations [JT].

Vortex equations over K\"ahler manifolds have been investigated  by Bradlow
[B1], [B2] and by Garcia-Prada [G1], [G2]: Given a pair $({\Cal E},\varphi)$
consisting of a holomorphic vector bundle with a section, the vortex
equations ask
for a Hermitian metric $h$ in ${\Cal E}$ with prescribed mean curvature: more
precisely, the equations---which depend on a real parameter
$\tau$---are
$$i\Lambda F_h=\frac{1}{2}(\tau\id_{\cal E}-\varphi\otimes\varphi^*).$$
A solution exists if and only if the pair $({\Cal E},\varphi)$ satisfies a
certain stability condition ($\tau$-stability), and the moduli space of
vortices can be identified with the moduli space of $\tau$-stable pairs. A GIT
construction of the latter space has been given by Thaddeus [T] and
Bertram [B] if the base manifold is a projective curve, and by Huybrechts
and Lehn
[HL1], [HL2] in the case of a projective variety. Other constructions have been
given by Bradlow and Daskalopoulus [BD1], [BD2] in the case of a Riemann
surface,
and by Garcia-Prada for compact K\"ahler manifolds [G2]. In this connection
also
[BD2] is relevant.
In this note we prove the following result:
\begin{th} Let $(X,g)$ be a K\"ahler surface of total scalar curvature
$\sigma_g$, and let $\Sigma$ be the canonical $\Spin^c$-structure with
associated Chern connection $c$. Fix a unitary vector bundle $E$ of rank $r$
over $X$, and define $\mu_g(\Sigma^+\otimes E):=\frac{\deg_g(E)}{r} +\sigma_g$.

Then for $\mu_g<0$, the space of solutions of the coupled
Seiberg-Witten equation is isomorphic to the moduli space of stable pairs of
topological type
$E$, with parameter $\sigma_g$.
\end{th}

If the constant $\mu_g(\Sigma^+\otimes E)$ is positive, one simply replaces the
bundle
$E$ with $E^{\vee}\otimes K_X$, where $K_X$ denotes the canonical line bundle
of
$X$ (cf. Lemma 3.1).

Note that the above theorem gives a complex geometric interpretation of the
moduli space of solutions of the coupled Seiberg-Witten equation associated
to \underbar{all} $\Spin^c$-structures on $X$: The change of the
$\Spin^c$-structure is equivalent to tensoring $E$ with a line bundle.

Notice also
that in the special case
$r=1$ one recovers Witten's result identifying the space of irreducible
monopoles
on a K\"ahler surface with the set of all divisors associated to
line bundles of a fixed topological type; the stability condition which he
mentions is the rank-1 version of the stable pair-condition.

Having established this correspondence, we describe some of the basic
properties of the moduli spaces, and give a first application: We show that
minimal surfaces of general type cannot be diffeomorphic to rational ones. This
provides a short proof of one of the essential steps in Friedman and Qin's
proof
of the Van de Ven conjecture [FQ]. More detailed investigations and
applications
will appear in a later paper.

We like to thank A. Van de Ven for very helpful questions and remarks.

\section{$\Spin^c$-structures and almost canonical classes}

The complex spinor group is defined as $\Spin^c:=\Spin\times_{\Z_2}S^1$, and
there are two non-split  exact sequences
$$\begin{array}{c}
1\map S^1\map\Spin^c\map\SO\map 1\\
1\map\Spin\map\Spin^c\map\ S^1\map 1\end{array}$$

In dimension 4, $\Spin^c(4)$ can be identified with the subgroup of
$\U(2)\times\U(2)$ consisting of pairs of unitary matrices with the same
determinant, and one has two commutative diagrams:

$$\begin{array}{ccllclll}
           &    &     1             &    &    1     &        &         &    \\

           &    &\downarrow         &    &\downarrow&        &          &    \\
          1&\map&  \Z_2             &\map&\Spin(4)  &\map    &\SO(4)\ \ \
&\map
1\\
           &    &\downarrow         &    &\downarrow&        & \ \
 \parallel&
\\
          1&\map&S^1                &\map&\Spin^c(4) &\map    &\SO(4) \ \ \
&\map 1\\
           &    &\downarrow(\cdot)^2&
           &{\scriptstyle{\det}}\downarrow\ \ \ \
&\ \nwarrow{\scriptstyle\Delta}&\ \ \uparrow&
\\
           &    &S^1
           &=   & S^1               &\longleftarrow&\U(2)\ \ \  &  \\
           &    &\downarrow         &     &\downarrow&{\ }^{\det}&&\\
           &    &1                  &    &1         &&&
\end{array} \eqno{(1)}$$
where $\Delta:\U(2)\map\Spin^c(4)\subset\U(2)\times\U(2)$ acts by
$a\map\left(\left(\matrix{\id&0\cr 0&\det a\cr}\right),a\right)$, and
$$
\matrix{
 & & &    &          &    &          &    &   1\ \ \     & &\cr
 & & &    &          &    &          &    &\downarrow\ \ \  & &\cr
 & & &    & 1        &    & 1        &    & \Z_2\ \ \       & &\cr
 & & &    &\downarrow&    &\downarrow&    &\downarrow\ \ \  & &\cr
 & &1&\rightarrow&S^1       &\rightarrow&\Spin^c(4) &\rightarrow&\SO(4)
&\rightarrow &1\cr
 & & &    &\downarrow&    &\downarrow&    &\ \
\  \ \ \ \ \downarrow{\scriptstyle(\lambda^+,\lambda^-)}& &\cr
 & &1&\rightarrow&S^1\times S^1
 &\rightarrow&\U(2)\times\U(2)&\stackrel{\rm
ad}{\rightarrow}&\SO(3)\times\SO(3)\ &\rightarrow&1\cr
 & & &    &\downarrow&    &\downarrow&    &\downarrow\ \ & &\cr
1&\rightarrow&\Z_2&\rightarrow&S^1&\stackrel{(\cdot)^2}{\rightarrow}&S^1&
\rightarrow&1\ \ &\cr
& & &    &\downarrow&    &\downarrow&    & \ \ \ \ & &\cr
& & &    & 1        &    & 1        &    &       & &\cr}$$
where $\lambda^{\pm}:\SO(4)\map\SO(3)$ are induced by the
two projections of $\Spin(4)=\SU(2)^+\times\SU(2)^-$ [HH]. $\lambda^{\pm}$ can
be also be seen as the representations of $\SO(4)$ in
$\Lambda^2_{\pm}(\R^4)\simeq\R^3$ induced by the canonical representation in
$\R^4$.

Let $X$ be a closed, oriented 4-manifold. Given any principal $\SO(4)$-bundle
$P$ over $X$, we denote by  $P^{\pm}$ the induced principal $\SO(3)$-bundles.
If
$\hat{P}$ is a $\Spin^c(4)$-bundle, we let $\Sigma^{\pm}$ be the
associated $\U(2)$-vector bundles, and we set (via the vertical determinant-map
in (1)) $\det(\hat{P})=L$, so that $\det(\Sigma^{\pm})=L$.
\begin{lm}
Let $P$ be a principal $SO(4)$-bundle over $X$ with characteristic classes
$w_2(P)\in H^2(X,\Z_2)$, and $p_1(P), e(P)\in H^4(X,\Z)$. Then\hfill{\break}
i) $P$ lifts to a principal $\Spin^c(4)$-bundle $\hat{P}$ iff $w_2(P)$ lifts to
an integral cohomology class.\hfill{\break}
ii) Given a class $L\in H^2(X,\Z)$ with $w_2(P)\equiv\bar L$(mod 2), the
$\Spin^c(4)$-lifts $\hat{P}$ of $P$ with $\det\hat{P}=L$ are in 1-1
correspondence with the 2-torsion elements in $H^2(X,\Z)$.\hfill{\break}
iii) Let $\hat P$ be a $\Spin^c(4)$-principal bundle with $P\simeq\hat{P}/S^1$,
and let $L=\det\hat{P}$. Then the Chern classes of $\Sigma^{\pm}$ are:
$$\begin{array}{rl}c_1(\Sigma^{\pm})&=L\\
c_2(\Sigma^{\pm})&=\frac{1}{4}\left(L^2 -p_1(P)\mp 2e(P)\right)\end{array}$$
\end{lm}
\pf [HH] and the diagrams above.
\qed

 Consider now a Riemannian metric $g$ on $X$, and let $P$ be the associated
principal $\SO(4)$-bundle. In this case the real vector bundles associated to
$P^{\pm}$  via the standard representations are the bundles $\Lambda^2_{\pm}$
of (anti-) self-dual 2-forms on $X$.

The integral characteristic classes of
$P$ are given by $p_1(P)=3\sigma$ and $e(P)=e$, where $\sigma$ and $e$
denote the
signature and the Euler number of the oriented manifold $X$. Furthermore,
$w_2(P)$ always lifts to an integral class, the lifts are precisely the
characteristic elements in $H^2(X,\Z)$, i.e. the classes $L$ with $x^2\equiv
x\cdot L$ for every $x\in H^2(X,\Z)$ [HH].

Let  $T_X$ be the tangent bundle of $X$, and denote by $\Lambda^p$ the bundle
of $p$-forms on $X$. The choice of a $\Spin^c(4)$-lift $\hat{P}$ of $P$ with
associated $\U(2)$-vector bundles $\Sigma^{\pm}$ defines a vector bundle
isomorphism
$\gamma:\Lambda^1\otimes\C\map\Hom_{\C}(\Sigma^+,\Sigma^-)$. There is also a
$\C$-linear isomorphism $(\cdot)^{\#}:\Hom_{\C}(\Sigma^+,\Sigma^-)\map
\Hom_{\C}(\Sigma^-,\Sigma^+)$ which satisfies the identity:
$$\gamma(u)^{\#}\gamma(v)+\gamma(v)^{\#}\gamma(u)=2g^{\C}(u,v)\id_{\Sigma^+},$$
and $\gamma(u)^{\#}=\gamma(u)^*=g(u,u)\gamma(u)^{-1}$ for real non-vanishing
cotangent vectors $u$.

It is convenient to extend the homomorphisms $\gamma(u)$ to endomorphisms of
the direct sum $\Sigma:=\Sigma^+\oplus\Sigma^-$. Putting
$\gamma(u)|_{\Sigma^-}:=-\gamma(u)^{\#}$, we obtain a vector-bundle
homomorphism $\gamma:\Lambda^1\otimes{\C}\map\End_0(\Sigma)$, which maps the
bundle $\Lambda^1$ of real 1-forms into the bundle of trace-free
skew-Hermitian endomorphisms of $\Sigma$. With this convention, we get:
$$\gamma(u)\circ\gamma(v)+\gamma(v)\circ\gamma(u)=-2g^{\C}(u,v)\id_{\Sigma}.$$

Consider the induced homomorphism
$$\Gamma:\Lambda^2\otimes\C \map\End_0(\Sigma)$$
defined on decomposable elements by
$$\Gamma(u\wedge v):=\frac{1}{2}[\gamma(u),\gamma(v)].$$
The restriction $\Gamma|_{\Lambda^2}$ identifies the bundle $\Lambda^2$
with the bundle $\ad_0(\hat{P})\simeq\ad(P)$ of skew-symmetric endomorphisms of
the tangent bundle of $X$.

$\Lambda^2$ splits as an orthogonal sum
$\Lambda^2=\Lambda^2_+\oplus\Lambda^2_-$ and $\Gamma$ maps
the bundle $\Lambda^2_{\pm}\otimes\C$ (respectively $\Lambda^2_{\pm}$)
isomorphically onto the bundle
$\End_0(\Sigma^{\pm})\subset\End(\Sigma)$
($su(\Sigma^{\pm})\subset su(\Sigma)$) of trace-free (trace free
skew-Hermitian)
endomorphisms of
$\Sigma^{\pm}$.

We give now an explicit description of the two spinor bundles $\Sigma^{\pm}$
and of the map $\Gamma$ in the case of a $\Spin^c(4)$-structure coming from an
almost Hermitian structure.
\begin{dt}
A characteristic element $K\in H^2(X,\Z)$ is an almost canonical class if
$K^2=3\sigma+2e$.
\end{dt}
Such classes exist on $X$ if and only if $X$ admits an almost complex
structure.
More precisely:
\begin{pr} ({\rm Wu}) $K\in H^2(X,\Z)$ is an almost canonical class if and only
if there exists an almost complex structure $J$ on $X$ which is compatible with
the orientation, such that $K=c_1(\Lambda^{10}_J)$.
\end{pr}
\pf [HH]
\qed

Here we denote, as usual, by $\Lambda^{pq}_J$ the bundle of $(p,q)$-forms
defined by the almost complex structure $J$.

Notice that any almost complex structure $J$ on $X$ can be deformed into a
$g$-orthogonal one, and that $J$ is $g$-orthogonal iff $g$ is $J$-Hermitian.
The choice of a $g$-orthogonal almost complex structure $J$ on $X$ corresponds
to to a reduction of the $\SO(4)$-bundle $P$ of $X$ to a $U(2)$-bundle via the
inclusion $\U(2)\subset\SO(4)$; since the inclusion factors through the
embedding $\Delta:\U(2)\map\Spin^c(4)$ (see diagram (1)), this
reduction defines a unique $\Spin^c(4)$-bundle $\hat{P_J}$ over $X$. By
construction we have $\hat{P}_J/{S^1}\simeq P$, and $\det\hat{P}_J=-K$.
\begin{pr} Let $J$ be a $g$-orthogonal almost complex structure on $X$,
compatible with the orientation.\hfill{\break}
i) The spinor bundles $\Sigma^{\pm}_J$ of $\hat{P}_J$ are:
$$\Sigma^+_J\simeq\Lambda^{00}\oplus\Lambda^{02}_J,\ \
\Sigma^-_J\simeq\Lambda^{01}_J.$$
ii) The map $\Gamma:\Lambda_+^2\otimes\C\map\End_0(\Sigma^+_J)$ is given by
$$\Lambda^{20}_J\oplus\Lambda^{02}_J\oplus\Lambda^{00}\omega_g
\ni(\lambda^{20},\lambda^{02},\omega_g)\stackrel{\Gamma}{\longmapsto}
2\left[\matrix{-i&-*(\lambda^{20}\wedge\cdot)\cr
\lambda^{02}\wedge\cdot&i\cr}\right]\in\End_0(\Lambda^{00}\oplus\Lambda^{02}).$$
\end{pr}
\pf i) $c_1(\Sigma^+_J)=c_1(\Sigma^-_J)=-K$,
$c_2(\Sigma^+_J)=\frac{1}{4}[K^2-3\sigma-2e]$,
$c_2(\Sigma^-_J)=\frac{1}{4}[K^2-3\sigma+2e]=c_2(\Sigma^+)+e$,
and $\U(2)$-bundles on a 4-manifold are classified by their Chern
classes.\hfill{\break}
ii) With respect to a suitable choice of the isomorphisms i), the Clifford map
$\gamma$ acts by
$$\gamma(u)(\varphi+\alpha)=\sqrt{2}\left(\varphi u^{01}-i\Lambda_g
u^{10}\wedge\alpha\right),$$
$$\gamma(u)^{\#}(\theta)=\sqrt{2}\left(i\Lambda_g(u^{10}\wedge\theta)-u^{01}
\wedge
\theta\right),
\eqno{(3)}$$
where $\Lambda_g:\Lambda^{pq}_J\map\Lambda^{p-1,q-1}_J$ is the adjoint of
the map
$\cdot\wedge\omega_g$ [H1].
\qed

\section{The coupled Seiberg-Witten equations}

Let $P$ be the principal $\SO(4)$-bundle associated with the tangent bundle
of the oriented, closed Riemannian 4-manifold $(X,g)$, and fix a $\Spin^c(4)$
structure $\hat{P}$ over $P$ with $L=\det(\hat{P})$. The choice of a
$\Spin^c(4)$-connection  in $\hat{P}$ projecting onto the Levi-Civita
connection in $P$ is equivalent to the choice of a connection $b$ in the
unitary line bundle $L$ [H1]. We denote by $B(b)$ the $\Spin^c(4)$-connection
in
$\hat{P}$ corresponding to $b$, and also the induced
connection in the vector bundle $\Sigma=\Sigma^+\oplus\Sigma^-$.
The curvature $F_{B(b)}$ of the connection $B(b)$ in $\Sigma$ has the form
$$F_{B(b)}=\frac{1}{2}F_b\id_{\Sigma}+F_g=
\left[\matrix{\frac{1}{2}F_b\id_{\Sigma^+}+F_g^+&0\cr
0&\frac{1}{2}F_b\id_{\Sigma^-}+F_g^-\cr}\right],$$
where $F_g$, and $F_g^{\pm}$ denote the Riemannian curvature operator, and its
components with respect to the splitting $\ad(P)=\Lambda^2_+\oplus\Lambda^2_-$.

Let now $E$ be an arbitrary Hermitian bundle of rank $r$ over $X$, and $A$ a
connection in $E$. We denote by $A_b$ the induced connection in the tensor
product $\Sigma\otimes E$, and by $\Dr_{A,b}:A^0(\Sigma\otimes E)\map
A^0(\Sigma\otimes E)$ the associated Dirac operator. $\Dr_{A,b}$ is defined as
the composition:
$$A^0(\Sigma\otimes E)\stackrel{\nabla_{A_b}}{\map}A^1(\Sigma\otimes
E)\stackrel{m}{\map} A^0(\Sigma\otimes E)$$
where $m$ is the Clifford multiplication $m(u,\sigma\otimes
e):=\gamma(u)(\sigma)\otimes e$. $\Dr_{A,b}$ is an elliptic, self-adjoint
operator and its square $\Dr_{A,b}^2$ is related to the usual Laplacian
$\nabla_{A_b}^*\nabla_{A_b}$ by the Weitzenb\"ock formula
$$\Dr_{A,b}^2=\nabla_{A_b}^*\nabla_{A_b}+\Gamma(F_{A_b}).$$
Here $\Gamma(F_{A_b})\in A^0(\End(\Sigma\otimes E))$ is the Hermitian
endomorphism defined as the composition
$$A^0(\Sigma\otimes E)\textmap{F_{A_b}}A^0(\Lambda^2\otimes\Sigma\otimes E)
\textmap{\Gamma} A^0(\End_0(\Sigma)\otimes\Sigma\otimes
E)\textmap{ev}A^0(\Sigma\otimes E).$$
We set $F_{A,b}:=F_A+\frac{1}{2}F_b\id_E\in A^0(\Lambda^2\otimes\End(E))$.

\begin{pr}
Let $s$ be the scalar curvature  of the Riemannian 4-manifold $(X,g)$.  Fix a
$\Spin^c(4)$-structure on $X$ and choose connections $b$ and $A$ in $L$ and $E$
respectively. Then
$$\Dr_{A,b}^2=\nabla_{A_b}^*\nabla_{A_b}+\Gamma(F_{A,b})+\frac{s}{4}
\id_{\Sigma\otimes E}.$$
\end{pr}
\pf Since $\Gamma(F_g)=\frac{s}{4}\id_{\Sigma}$ [H1], and
$F_{A_b}=F_{B(b)}\otimes
\id_E+\id_{\Sigma}\otimes
F_A=\frac{1}{2}F_b\id_{\Sigma}\otimes\id_E+F_g\otimes\id_E+\id_{\Sigma}\otimes
F_A=\id_{\Sigma}\otimes(F_A+\frac{1}{2}F_b\id_E)+F_g\id_E$, we find
$\Gamma(F_{A_b})=\Gamma(F_{A,b})+\frac{s}{4}\id_{\Sigma\otimes E}$.
\qed

\begin{re} One has a Bochner-type result for spinors $\Psi$ on which\linebreak
\hbox{$\Gamma(F_{A,b})+\frac{s}{4}\id_{\Sigma\otimes E}$}  is positive: Such a
spinor is harmonic if and only if it is parallel [H1].
\end{re}

Let $(\ ,\ )$  be the pointwise inner product on $\Sigma\otimes E$, $|\ |$
the associated pointwise norm, and $\parallel\ \parallel$ the corresponding
$L^2$-norm. For a spinor
$\Psi\in A^0(\Sigma^{\pm}\otimes E)$ we define
$(\Psi\bar\Psi)_0\in A^0(\End_0(\Sigma^{\pm}\otimes E))$ as the image of the
Hermitian endomorphism $\Psi\otimes\bar\Psi\in A^0(\End(\Sigma^{\pm}\otimes
E))$
under the projection $\End(\Sigma^{\pm}\otimes
E)\map\End_0(\Sigma^{\pm})\otimes\End(E)$.
\begin{co}
With the notations above, we have
$$(\Dr_{A,b}^2\Psi,\Psi)=(\nabla_{A_b}^*\nabla_{A_b}\Psi,\Psi)+(\Gamma(F_{A,
b}^+),
(\Psi_+\bar\Psi_+)_0)+(\Gamma(F_{A,b}^-),
(\Psi_-\bar\Psi_-)_0)+\frac{s}{4}|\Psi|^2,$$
where ($F_{A,b}^-$) $F_{A,b}^{+}$ is the (anti-)self-dual component of
$F_{A,b}$.
\end{co}
\pf Indeed, since $\Gamma(F_{A,b}^{\pm})$ vanishes on $\Sigma^{\mp}$ and is
trace
free with respect to
$\Sigma^{\pm}$, the inner product
$(\Gamma(F_{A,b}),(\Psi\bar\Psi))$ in the Weitzenb\"ock formula simplifies for
a
spinor
$\Psi\in A^0(\Sigma^{\pm}\otimes E) $:
$$(\Gamma(F_{A,b}),(\Psi\bar\Psi))=(\Gamma(F_{A,b}^{\pm}),(\Psi\bar\Psi)_0)$$
\qed

For a positive spinor $\Psi\in A^0(E\otimes\Sigma^+)$, the following
important identity follows immediately:
$$(\Dr_{A,b}^2\Psi,\Psi)+\frac{1}{2}|\Gamma(F_{A,b}^+)-(\Psi\bar\Psi)_0|^2=
(\nabla_{A_b}^*\nabla_{A_b}\Psi,\Psi)+
\frac{1}{2}|F_{A,b}^+|^2+\frac{1}{2}|(\Psi\bar\Psi)_0|^2+\frac{s}{4}|\Psi|^2
\eqno{(4)}$$

If we integrate both sides of (4) over $X$, we get:
\begin{pr} Let $(X,g)$ be an oriented, closed Riemannian 4-manifold with scalar
curvature $s$, $E$ a Hermitian bundle over $X$. Choose a $\Spin^c(4)$-structure
on $X$ and a  connection $b$ in the determinant line bundle
$L=\det(\Sigma^+)=\det(\Sigma^-)$. Let $A$ be a connection in $E$. For any
$\Psi\in A^0(\Sigma^+\otimes E)$ we have:
$$\parallel\Dr_{A,b}\Psi\parallel^2+
\frac{1}{2}\parallel\Gamma(F_{A,b}^+)-(\Psi\bar\Psi)_0\parallel^2=$$ $$=
\parallel\nabla_{A_b}\Psi\parallel^2+
\frac{1}{2}\parallel
F_{A,b}^+\parallel^2+\frac{1}{2}\parallel(\Psi\bar\Psi)_0\parallel^2+
\frac{1}{4}\int\limits_X s|\Psi|^2.$$
\end{pr}

We introduce now our coupled variant of the Seiberg-Witten equations. The
unknown is a pair $(A,\Psi)$ consisting of a connection in the Hermitian bundle
$E$ and a section $\Psi\in A^0(\Sigma^+\otimes E)$. The equations ask for the
vanishing of the left-hand side in the above formula.
$$\left\{\begin{array}{ccc}\Dr_{A,b}\Psi&=&0\\
\Gamma(F_{A,b}^+)&=&(\Psi\bar\Psi)_0
\end{array}\right.\eqno{(SW)}$$
Proposition 2.4 and the inequality
$|(\Psi\bar\Psi)_0|^2\geq\frac{1}{2}|\Psi|^4$
give immediately:
\begin{re} If the scalar curvature $s$ is nonnegative on $X$,  then the only
solutions of the equations are the pairs $(A,0)$, with $F_{A,b}^+=0$.
\end{re}

If $L$ is the square of a line bundle $L^{\frac{1}{2}}$, and if we choose a
connection $b^{\frac{1}{2}}$ in $L^{\frac{1}{2}}$ with square $b$, then
$F_{A,b}$ is simply the curvature of the connection $A_{b^{\frac{1}{2}}}$ in
$E\otimes L^{\frac{1}{2}}$. The solution of the coupled Seiberg-Witten
equations
on a manifold with $s\geq 0$ are in this case just $\U(r)$-instantons on
$E\otimes L^{\frac{1}{2}}$.

In the case of a K\"ahler surface $(X,g)$, the coupled Seiberg-Witten equation
can be reformulated in terms of complex geometry. The point is that if we
consider the canonical $\Spin^c(4)$-structure associated to the K\"ahler
structure, the Dirac operator has a very simple form [H1]. The determinant of
this $\Spin^c(4)$-structure is the anti-canonical bundle $K_X^{\vee}$ of the
surface, which comes with a holomorphic structure and a natural metric
inherited from the holomorphic tangent bundle.

Let $c$ be the Chern connection in $K_X^{\vee}$. With this choice, the induced
connection $B(c)$ in
$\Sigma=\Lambda^{00}\oplus\Lambda^{02}\oplus\Lambda^{01}$ coincides with the
connection defined by the Levi-Civita connection. Recall that on a K\"ahler
manifold, the almost complex structure is parallel with respect to the
Levi-Civita connection, so that the splitting of the exterior algebra
$\bigoplus\limits_{p}\Lambda^p\otimes\C$ becomes parallel, too.
\begin{pr}
Let $(X,g)$ be a K\"ahler surface with Chern connection $c$ in $K_X^{\vee}$.
Choose a connection $A$ in a Hermitian vector bundle
 $E$ over $X$ and a section $\Psi=\varphi+\alpha\in
A^0(E)+A^0(\Lambda^{02}\otimes E)$.

The pair $(A,\Psi)$ satisfies the Seiberg-Witten equations iff the following
identities hold:
$$
\begin{array}{lll}
F_{A,c}^{20}&=&-\frac{1}{2}\varphi\otimes\bar\alpha\\
F_{A,c}^{02}&=&\frac{1}{2}\alpha\otimes\bar\varphi\\
i\Lambda_g F_{A,c}&=&-\frac{1}{2}\left(\varphi\otimes\bar\varphi-
*(\alpha\otimes\bar\alpha)\right)\\
\bar\partial_A\varphi&=&i\Lambda_g\partial_A\alpha\end{array}$$
\end{pr}
\pf The Dirac operator is in this case
$\Dr_{A,c}=\sqrt{2}(\bar\partial_A-i\Lambda_g\partial_A)$, and the endomorphism
$(\Psi\bar\Psi)_0$ has the form:
$$\left[\matrix{\frac{1}{2}(\varphi\otimes\bar\varphi-
*(\alpha\otimes\bar\alpha))&*(\varphi\otimes\bar\alpha\wedge\cdot)\cr
\alpha\otimes\bar\varphi&-\frac{1}{2}(\varphi\otimes\bar\varphi-
*(\alpha\otimes\bar\alpha))\cr}\right].$$

Since $\Gamma(F_{A,c}^+)=\Gamma(F_{A,c}^{20}+F_{A,c}^{02}+\frac{1}{2}\Lambda_g
F_{A,c}\cdot\omega_g)$ equals
$$2\left[\matrix{-\frac{i}{2}\Lambda_gF_{A,c}&-*(F_{A,c}^{20}\wedge\cdot)\cr
F_{A,c}^{20}\wedge\cdot&\frac{i}{2}\Lambda_gF_{A,c}\cr}\right],$$
the equivalence of the two systems of equations follows.
\qed

\section{Monopoles on K\"ahler surfaces and the generalized vortex equation}

Let $(X,g)$ be a K\"ahler surface with canonical $\Spin^c(4)$-structure, and
Chern connection $c$ in the anti-canonical bundle $K_X^{\vee}$.

We fix a unitary vector bundle $E$ of rank $r$ over $X$, and define
$J(E):=\deg_g(\Sigma^+\otimes E)$, i.e.
$J(E)=2r(\mu_g(E)-\frac{1}{2}\mu_g(K_X))$, where $\mu_g$ denotes the slope with
respect to $\omega_g$.

Every spinor $\Psi\in A^0(\Sigma^+\otimes E)$ has the form
$\Psi=\varphi+\alpha$ with $\varphi\in A^0(E)$ and $\alpha\in
A^{0}(\Lambda^{02}\otimes E)$.

We have seen that the coupled Seiberg-Witten equations are equivalent to the
system:
$$\left\{
\begin{array}{lll}
F_{A,c}^{20}&=&-\frac{1}{2}\varphi\otimes\bar\alpha\\
F_{A,c}^{02}&=&\frac{1}{2}\alpha\otimes\bar\varphi\\
i\Lambda_g F_{A,c}&=&-\frac{1}{2}\left(\varphi\otimes\bar\varphi-
*(\alpha\otimes\bar\alpha)\right)\\
\bar\partial_A\varphi&=&i\Lambda_g\partial_A\alpha\end{array}\right.
\eqno{(SW^*)}$$

\begin{lm} \hfill{\break}
A. Suppose $J(E)<0$: \hfill{\break}
A pair $(A,\varphi+\alpha)$ is a solution of the system $(SW^*)$ if and only
if \hfill{\break}
i) $F_A^{20}=F_A^{02}=0$\hfill{\break}
ii) $\alpha=0$, $\bar\partial_A\varphi=0$ \hfill{\break}
iii) $i\Lambda_g F_A+\frac{1}{2}\varphi\otimes\bar\varphi+\frac{1}{2}s\id_E=0$.
\hfill{\break}
B. Suppose $J(E)>0$, and put $a:=\bar\alpha\in A^{20}(\bar E)=A^0(E^{\vee}
\otimes K_X)$:\hfill{\break}
A pair $(A,\varphi+\bar a)$ is a
solution of the system
$(SW^*)$ if and only if\hfill{\break}
i) $F_A^{20}=F_A^{02}=0$\hfill{\break}
ii) $\varphi=0$, $\bar\partial_A a=0$ \hfill{\break}
iii) $i\Lambda_g F_A-\frac{1}{2}*(a\otimes\bar a)+\frac{1}{2}s\id_E=0$.

\end{lm}
\pf (cf. [W]) The splitting $\Sigma^+\otimes E=\Lambda^{00}\otimes
E\oplus\Lambda^{02}\otimes E$ is parallel with respect to $\nabla_{A,c}$, so
that, by Proposition 2.4
$$\parallel\Dr_{A,c}\Psi\parallel^2+
\frac{1}{2}\parallel\Gamma(F_{A,c}^+)-(\Psi\bar\Psi)_0\parallel^2=$$ $$=
\parallel\nabla_{A_c}\varphi\parallel^2+
\parallel\nabla_{A_c}\alpha\parallel^2+
\frac{1}{2}\parallel
F_{A,c}^+\parallel^2+\frac{1}{2}\parallel(\Psi\bar\Psi)_0\parallel^2+
\frac{1}{4}\int\limits_X s(|\varphi|^2+|\alpha|^2).$$
The right-hand side is invariant under the transformation
$(A,\varphi,\alpha)\longmapsto (A,\varphi,-\alpha)$, hence any solution
$(A,\varphi+\alpha)$ must have $F_A^{20}=F_A^{02}=0$ and
$\varphi\otimes\bar\alpha=\alpha\otimes\bar\varphi=0$; the latter implies
obviously $\alpha=0$ or $\varphi=0$. Integrating the trace of the equation
$i\Lambda F_{A,c}=-\frac{1}{2}\left(\varphi\otimes\bar\varphi-
*(\alpha\otimes\bar\alpha)\right)$, we find:
$$J(E)=c_1(\Sigma^+\otimes
E)\cup[\omega_g]=(2c_1(E)-rc_1(K_X))\cup[\omega_g]=$$ $$=
2\int\limits_X\frac{i}{2\pi}\tr(F_{A,c})\wedge\omega_g=
\frac{1}{4\pi}\int\limits_X\tr(i\Lambda
F_{A,c})\omega_g^2=\frac{1}{8\pi}\int\limits_X\tr(-\varphi\otimes\bar\varphi)
+*(\alpha\otimes\bar\alpha))\omega_g^2$$
This equation shows that we must have $\alpha=0$, if $J(E)<0$, and
$\varphi=0$, if $J(E)>0$.
Notice
that, replacing   $E$ by $E^{\vee}\otimes K_X$, the second case reduces to the
first one.

The assertion follows now immediately from the identity
$i\Lambda_g F_c=s$.

\qed

Notice that the last equation
$$i\Lambda_g F_A+\frac{1}{2}\varphi\otimes\bar\varphi+\frac{1}{2}s\id_E=0$$
has the form of a generalized vortex equation as studied by Bradlow [B1], [B2]
and by Garcia-Prada [G2]; it is precisely the vortex equation with constant
$\tau=-{s}$, if $(X,g)$ has constant scalar curvature.

Let $s_m$ be the mean scalar curvature defined by
$\int\limits_Xs\omega_g^2=s_m\int\limits_X\omega^2=2s_m\Vol_g(X)$.

We are going to prove that the system
$$\left\{\begin{array}{cl}\bar\partial_A^2&=0\\ \bar\partial_A\varphi&=0\\
i\Lambda_g
F_A+\frac{1}{2}\varphi\otimes\bar\varphi+\frac{1}{2}s\id_E&=0
\end{array}\right.$$
for the pair $(A,\varphi)$ consisting of a unitary connection in $E$, and a
section in $E$, is always equivalent to the vortex system with parameter
$\tau=-s_m$, i.e. to the system obtained by replacing the third equation with
$$i\Lambda_g F_A+\frac{1}{2}\varphi\otimes\bar\varphi+\frac{1}{2}s_m\id_E=0.$$
"Equivalent" means here that the corresponding moduli spaces of solutions are
naturally isomorphic.

Let generally $t$ be a smooth real function on $X$ with mean value $t_m$, and
consider the following system of equations:
$$\left\{\begin{array}{cl}\bar\partial_A^2&=0\\ \bar\partial_A\varphi&=0\\
i\Lambda_g F_A+\frac{1}{2}\varphi\otimes\bar\varphi-\frac{1}{2}t\id_E&=0
\end{array}\right.\eqno(V_t)$$
$(V_t)$ is defined on the space $\AA(E)\times A^0(E)$, where $\AA(E)$ is the
space of unitary connections in $E$. The product $\AA(E)\times A^0(E)$
(completed with respect to sufficiently large  Sobolev indices) carries a
natural
$L^2$ K\"ahler metric $\tilde g$ and a natural right action of the gauge group
$U(E)$: $(A,\varphi)^f:=(A^f,f^{-1}\varphi)$, where
$d_{A^f}:=f^{-1}\circ d_A\circ f$.

For every real function $t$ let
$$m_t:\AA(E)\times A^0(E)\map A^0(\ad(E))$$
be the map given by $m_t:=\Lambda_g
F_A-\frac{i}{2}\varphi\otimes\bar\varphi+\frac{i}{2}t\id_E.$
\begin{pr}
$m_t$ is a moment map for the action of $U(E)$ on \linebreak $\AA(E)\times
A^0(E)$.
\end{pr}
\pf  Let $a^{\#}$ be the vector field  on $\AA(E)\times A^0(E)$ associated
with the infinitesimal transformation $a\in A^0(\ad(E))={\rm Lie }(U(E))$, and
define the real function $m^a_t:\AA(E)\times A^0(E)\map\R$ by:
$$ m^a_t(x)=\langle m_t(x),a\rangle_{L^2} .$$
 We
have to show that
$m_t$ satisfies the identities:
$$\iota_{a^{\#}}\omega_{\tilde g}=dm_t^a  \ ,\ \ \   m_t^a\circ
f=m^{\ad_f(a)}\ \ \ {\rm for\ all}\
\ a\in A^0(\ad(E)), \ \ f\in U(E).$$

It is well known that, in general, a moment map for a group action in a
symplectic
manifold is well defined up to a constant central element in the Lie algebra
of the group. In our case, the center of the Lie algebra
$A^0(\ad(E))$ of the gauge group is just $iA^0\id_E$, hence it suffices to show
that $m_0$ is a moment map. This has already been noticed by Garcia-Prada [G1],
[G2].
\qed

Note also that in our case every moment map has the form $m_t$ for some
function
$t$, which shows that from the point of view of symplectic geometry, the
natural
equations are the generalized vortex equations
$(V_t)$.

In order to show that Bradlow's main result [B2] also holds for the
generalized system $(V_t)$, we have to recall some definitions.

Let $\EE$ be a holomorphic vector bundle of topological type $E$, and let
$\varphi\in H^0(\EE)$ be a holomorphic section. The pair $(\EE,\varphi)$ is
$\lambda$-\underbar{stable} with respect to a constant $\lambda\in\R$ iff the
following conditions hold:\\
(1) $\mu_g(\EE)<\lambda$ and $\mu_g(\FF)<\lambda$ for all reflexive subsheaves
$\FF\subset\EE$ with $0<\rk(\FF)<r$.\\
(2) $\mu_g(\EE/\FF)>\lambda$ for all reflexive subsheaves
 $\FF\subset\EE$ with $0<\rk(\FF)<r$ and $\varphi\in H^0(\FF)$.

\begin{th} Let $(X,g)$ be a closed K\"ahler manifold, $t\in A^0$ a real
function, and $(\EE,\varphi)$ a holomorphic pair over $X$. Set
$\lambda:=\frac{1}{4\pi} t_m\Vol_g(X)$. $\EE$ admits a Hermitian metric $h$
such that the associated Chern connection $A_h$ satisfies the vortex equation
$$i\Lambda_g F_A+\frac{1}{2}\varphi\otimes\bar\varphi-\frac{1}{2}t\id_E=0$$
iff one of the following conditions holds:\\
(i) $(\EE,\varphi)$ is $\lambda$-stable\\
(ii) $\EE$ admits a splitting $\EE=\EE'\oplus\EE''$ with $\varphi\in
H^0(\EE')$ such that $(\EE',\varphi)$ is $\lambda$-stable, and $\EE''$ admits
a weak Hermitian-Einstein metric with factor $\frac{t}{2}$. In particular
$\EE''$ is polystable of slope $\lambda$.
\end{th}

\pf In the case of a constant function $t=\tau$, the theorem was proved by
Bradlow [B2], and his arguments work in the general context, too: The fact
that the existence of a solution of the vortex equation implies $(i)$ or
$(ii)$ follows by replacing the constant $\tau$ in [B2] everywhere with the
function $t$. The difficult part consists in  showing that every
$\lambda$-stable pair $(\EE,\varphi)$ admits a metric $h$ such that
$(A_h,\varphi)$ satisfies the vortex equation $(V_t)$. To this end let
$Met(E)$ be the space of Hermitian metrics in $E$, and fix a background
metric $k\in Met(E)$. Bradlow constructs a functional
$M_{\varphi,\tau}(\cdot,\cdot):Met(E)\times Met(E)\map\R$, which is convex in
the second argument, such that any critical point of
$M_{\varphi,\tau}(k,\cdot)$ is a solution of the vortex equation; the point
is then to find an absolute minimum of $M_{\varphi,\tau}(k,\cdot)$. The
existence of an absolute minimum follows from the following basic
$\CC^0$ estimate:
\begin{lm}{\rm (Bradlow)} Let $Met^p_2(E,B):=\{h=ke^a| a\in L^2_p(\End(E)),
a^{*k}=a, \parallel\mu_t(A_h,\varphi)\parallel_{L^p}\leq B\}$. If
$(\EE,\varphi)$ is $\frac{\tau}{4\pi}\Vol_g(X)$-stable, then there exist
positive constants $C_1$, $C_2$ such that
$$\sup|a|\leq C_1M_{\varphi,\tau}(k,ke^a)+ C_2,$$
for all $k$-Hermitian endomorphisms $a\in L^2_p(\End(E))$. Moreover, any
absolute minimum of $M_{\varphi,t}(k,\cdot)$ on $Met^p_2(E,B)$ is a critical
point of $M_{\varphi,t}(k,\cdot)$, and gives a solution of the vortex
equation $V_{\tau}$.
\end{lm}

Let now $t$ be a real function on $X$, and choose a solution $v$ of the
Laplace equation $i\Lambda_g\bar\partial\partial v=\frac{1}{2}(t-t_m)$. If we
make the substituion $h=h'e^v$, then $h$ solves the vortex equation $(V_t)$
 iff $h'$ is a solution of
$$i\Lambda_g
F_{h'}+\frac{1}{2}e^v\varphi\otimes\bar\varphi^{h'}-\frac{1}{2}t_m\id_E=0.$$
Define $\mu_{t_m,v}(h'):=i\Lambda_g
F_{h'}+\frac{1}{2}e^v\varphi\otimes\bar\varphi^{h'}-\frac{1}{2}t_m\id_E=0$, and
$$M_{\varphi,t_m,v}(k,h):=M_D(k,h)+\parallel
e^{\frac{v}{2}}\varphi\parallel^2_h-\parallel
e^{\frac{v}{2}}\varphi\parallel^2_k-t_m\int\limits_X\tr(\log(k^{-1}h)),$$
where $M_D$ is the Donaldson functional [D]. Then it is not difficult to show
that all arguments of Bradlow remain correct after replacing $\mu_{t_m}$ and
$M_{\varphi,t_m}$ with $\mu_{t_m,v}$ and $M_{\varphi,t_m,v}$ respectively.
Indeed, let $l$ be a positive bound from below for the map $e^v$. Then
$$\begin{array}{ll}M_{\varphi,t_m}(k,ke^{a+\log l})&\leq
M_D(k,ke^{a})+M_D(ke^{a},lke^{a})+\parallel
l\varphi\parallel^2_h-t_m\int\limits_X\tr\log(lk^{-1}h)\cr
&\leq M_{\varphi,t_m,v}(k,ke^{a})+\parallel
e^{\frac{v}{2}}\varphi\parallel^2_k+2\log l\deg_g(E)-rt_m\log l\Vol_g(X)\cr
&\leq  M_{\varphi,t_m,v}(k,ke^{a})+C'(k,\varphi,v,l).\end{array}$$
Similarly, we get constants $n>0$, $C''$ and an inequality
$$M_{\varphi,t_m,v}(k,ke^{a+\log n})\leq M_{\varphi,t_m}(k,ke^{a})+C'',$$
which shows that the basic $\CC^0$ estimate in the Lemma above holds for
$M_{\varphi,t_m,v}$ iff it holds for Bradlow's functional
$M_{\varphi,t_m}$.
\qed
\begin{re}
In the special case of a rank-1 bundle $E$, a much more elementary proof
based on [B1] is possible.
\end{re}

\section{Moduli spaces of monopoles, vortices, and stable pairs}

Let $(X,g)$ be a closed K\"ahler manifold of  arbitrary dimension, and fix a
unitary vector bundle $E$ of rank $r$ over $X$. We denote by $\bar{\AA}(E)$
the affine space of semiconnection  of type $(0,1)$ in $E$. the complex gauge
group $GL(E)$ acts on $\bar{\AA}(E)\times A^0(E)$ from the right by
$(\bar\partial_A,\varphi)^g:=(g^{-1}\circ\bar\partial_A\circ g,
g^{-1}\varphi)$;
this action becomes complex analytic after suitable Sobolev completions. We
denote
by
$\bar\SS(E)$ the set of pairs $(\bar\partial_A,\varphi)$ with trivial isotropy
group. Notice that $\varphi\ne 0$ when
$(\bar\partial_A,\varphi)\in\bar\SS(E)$,
and that  $\bar\SS(E)$ is an open subset of $\bar{\AA}(E)\times A^0(E)$, by
elliptic semi-continuity [K].

The action of $GL(E)$ on $\bar\SS(E)$ is free, by definition, and we denote
the Hilbert manifold $\qmod{\bar\SS(E)}{GL(E)}$ by $\bar{\BB}^s(E)$.
The map $p:\bar{\AA}(E)\times A^0(E)\map A^{02}(\End(E)\oplus A^{01}(E)$
defined by  $p(\bar\partial_A,\varphi)=(F_A^{02},\bar\partial_A\varphi)$ is
equivariant with respect to the natural actions of $GL(E)$, hence it gives
rise to a section $\hat p$ in the associated vector bundle
$\bar\SS(E)\times_{GL(E)}\left(A^{02}(\End(E)\oplus A^{01}(E)\right)$ over
$\bar{\BB}^s(E)$. We define the moduli space of \underbar{simple pairs} of type
$E$ to be the zero-locus $Z(\hat{p})$ of this section. $Z(\hat{p})$ can be
identified with the set of isomorphism classes consisting of a holomorphic
bundle $\EE$ of differentiable type $E$, and a holomorphic section
$\varphi\ne 0$, such that the kernel of the evaluation map
$ev(\varphi):H^0(\End(\EE))\map H^0(\EE)$ is trivial.

 In a similar way we define the moduli space
$\VV^g_t$ of gauge-equivalence classes of irreducible  solutions of the
generalized vortex equation $V_t$:

Let $B^+$ denote as usual the subbundle
$\left((\Lambda^{02}+\Lambda^{20})\cap \Lambda^2\right)\oplus
\Lambda^0\omega$ of
the bundle $\Lambda^2$ of real 2-forms on $X$. We denote by ${\Cal D}^*$ the
open subset of the product $\DD:=\AA(E)\times A^0(E)\simeq\bar{\AA}(E)\times
A^0(E)$ consisting of pairs with trivial isotropy group with respect to the
action of the gauge group $U(E)$. The quotient
$\BB^*(E):=\qmod{\DD^*(E)}{U(E)}$ comes with the structure of a real-analytic
manifold.

Let $v:\DD(E)\map A^0(B^+\otimes\ad(E))\oplus A^{01}(E)$ be the map given by:
$$v(A,\varphi)=(F^{20}+F^{02},m_t(A,\varphi)\omega_g\id_E,\bar\partial_A\varphi).$$
Again $v$ is $U(E)$-equivariant, and the moduli space
$\VV^g_t$  of $t$-vortices is defined to be the zero-locus $Z(\hat{v})$ of the
induced section $\hat v$  of \linebreak
\hbox{$\DD^*(E)\times_{U(E)}A^0(B^+\otimes\ad(E))\oplus A^{01}(E)$} over
$\BB^*(E)$.

Notice now that by Proposition 3.2,  the second component $v^2$ of $v$ is  a
moment map for the $U(E)$ action. It is easy to see that (at least
in a neighbourhood of $Z(v)\cap\DD^*$) it has the general
property of a moment map in the finite dimensional K\"ahler geometry: Its zero
locus $Z(v^2)$ is smooth and intersects every $GL(E)$ orbit along a $U(E)$
orbit,
and the intersection is transversal. This means that the natural map
$A\map\bar\partial_A$ defines in a neighbourhood of $Z(\hat{v})\cap\BB^*(E)$ an
open embedding $i:Z({\hat{v}^2})\map\bar\BB^s$ of smooth Hilbert manifolds.

Regard now $\VV^g_t$  as the subspace of $Z({\hat{v}^2})\subset\BB^*(E)$
defined by the equation $(\hat{v}^1,\hat{v}^3)=0$.  On the other hand, the
pullback of the equation
$\hat p=0$, cutting out the moduli space $Z(\hat{p})$ of simple holomorphic
pairs,
via  the open embedding $i$ is precisely the equation
$(\hat{v}^1,\hat{v}^3)=0$, cutting out $\VV^g_t$. We get therefore an open
embedding $i_0:\VV^g_t\map Z(\hat{p})$ of real analytic spaces induced by
$i$, and by Theorem 3.3 the image of $i_0$ consists of the set of
$\lambda$-stable pairs, with $\lambda:=\frac{1}{4\pi}t_m\Vol_g(X)$.

Finally we denote by $\MM_X^g(E,\lambda)\subset Z(\hat p)$ the open subspace of
$\lambda$-stable pairs, with the induced complex space-structure. Putting
everything together, we have:
\begin{th}
Let $(X,g)$ be a closed K\"ahler manifold, $t\in A^0$ a real function, and
$\lambda: =\frac{1}{4\pi}t_m\Vol_g(X)$. Fix a unitary vector bundle $E$ of rank
$r$ over $X$. There are natural real-analytic isomorphisms of moduli spaces
$$\VV^g_t(E)\simeq\VV^g_{t_m}(E)\simeq\MM_X^g(E,\lambda).$$
\end{th}

Let us come back now to the monopole equation $(SW^*)$ on a K\"ahler surface.
In
this case the function $t$ is the negative of the scalar curvature $s$, so that
the corresponding constant $\lambda$ becomes:
$$\lambda=\frac{-s_m}{4\pi}\Vol_g(X)=-\frac{1}{8\pi}\int\limits_Xs\omega^2=
-\frac{1}{8\pi}\int\limits_X(i\Lambda
F_c)\omega^2=-\frac{1}{4\pi}\int\limits_X i
F_c\wedge\omega=\frac{1}{2}\mu_g(K).$$

This yields our main result:
\begin{th}
Let $(X,g)$ be a K\"ahler surface with canonical $\Spin^c(4)$-structure, and
Chern connection $c$ in $K_X^{\vee}$. Fix a unitary vector bundle $E$ of
rank $r$
over $X$, and suppose $J(E)=\deg_g(\Sigma^+\otimes E)<0$. The moduli space of
solutions of the coupled Seiberg-Witten equations is isomorphic to the moduli
space $\MM_X^g(E,\frac{1}{2}\mu_g(K))$ of $\frac{1}{2}\mu_g(K)$-stable pairs of
topological type $E$.
\end{th}

At this point it is natural to study the properties of the moduli spaces
$\MM^g_X(E,\lambda)$. We do not want to go into details here, and we content
ourselves  by describing some of the basic steps.

The infinitesimal structure of the moduli space around a point
$[(A,\varphi)]$ is
given by a deformation complex
$(C_{\bar\partial_A,\varphi}^*, d_{A,\varphi}^*)$ which is the cone over the
evaluation map $ev^*$, $ev^q(\varphi):A^{0q}(\End(E))\map A^{0q}(E)$. More
precisely $C_{\bar\partial_A,\varphi}^q=A^{0q}(\End(E))\oplus A^{0,q-1}(E)$ and
the differential $d_{A,\varphi}^q$ is given by the matrix
$$d_{A,\varphi}^q=\left[\matrix{-\bar
D_A&0\cr ev(\varphi)&\bar\partial_A\cr}\right],$$
where $\bar\partial_A$ and $\bar D_A$ are the operators of the Dolbeault
complexes $(A^{0*}(E),\bar\partial_A)$ and $(A^{0*}\End(E),\bar D_A)$
respectively.

Associated to the morphism $ev^*(\varphi)$ is an exact sequence
$$\dots\map H^q(\End(\EE_A))\textmap{ev^q(\varphi)}H^q(\EE_A)\map
H_{\bar\partial_A,\varphi}^{q+1}\map H^{q+1}(\End(\EE_A))\map\dots $$
with finite dimensional vector spaces
$$H^q_{\bar\partial_A,\varphi}=
\ker(ev^q(\varphi))\oplus\coker(ev^{q-1}(\varphi)).$$

$H^0_{\bar\partial_A,\varphi}$ vanishes for a simple pair
$(\bar\partial_A,\varphi)$, and $H^1_{\bar\partial_A,\varphi}$ is the Zariski
tangent space of $\MM^g_X(E,\lambda)$ at $[\bar\partial_A,\varphi]$.

A Kuranishi type argument yields local models of the moduli space, which can be
locally described  as the zero loci of holomorphic map germs
$$K_{[\bar\partial_A,\varphi]}:H^1_{\bar\partial_A,\varphi}\map
H^2_{\bar\partial_A,\varphi}$$
at the origin.

One finds that $H^2_{\bar\partial_A,\varphi}=0$ is a sufficient smoothness
criterion in the point $[\bar\partial_A,\varphi]$ of the moduli space, and that
the expected dimension is \linebreak\hbox{$\chi(E)-\chi(\End(E))$}. The
necessary
arguments are very similar to the ones in [BD1], [BD2].

The moduli spaces $\MM^g(E,\lambda)$ will be quasi-projective varieties if the
underlying manifold $(X,g)$ is Hodge, i.e. if $X$ admits a projective embedding
such that a multiple of the K\"ahler class is a polarisation [G1].

A GIT construction for projective varieties of any dimension has been given in
[HL2]. The spaces $\MM^g_X(E,\lambda)$ vary with the
parameter
$\lambda$, and flip-phenomena occur just like in the case of curves [T].

\section{Applications}

The equations considered by Seiberg and Witten are associated to a
$\Spin^c(4)$-structure, and correspond to the case when (in our notations) the
unitary bundle
$E$ is the trivial line bundle. Alternatively, we can fix a $\Spin^c(4)$
structure $\s_0$ on $X$ , and  regard the Seiberg-Witten equations
corresponding
to the other $\Spin^c(4)$-structures as {\sl coupled} Seiberg-Witten
equations associated to $\s_0$ and to a unitary line bundle $E$. The
$\Spin^c(4)$-structure we fix will  always be the canonical structure
defined by a
K\"ahler metric. In the most interesting case of rank-1 bundles
$E$ over K\"ahler surfaces the central result is:
\begin{pr}
Let $(X,g)$ be a K\"ahler surface with canonical class $K$, and let $L$ be a
complex line bundle over $X$ with $L\equiv K$ (mod 2). Denote by
$\WW_X^g(L)$ the
moduli space of  solutions of the Seiberg-Witten equation for all
$\Spin^c(4)$-structures with determinant $L$. Then\hfill{\break}
i) If $\mu(L)<0$, $\WW_X^g(L)$ is isomorphic to the space of all linear systems
$|D|$, where $D$ is a divisor with $c_1(\OO_X(2D-K))=L$.\hfill{\break}
ii) If $\mu(L)>0$, $\WW_X^g(L)$ is isomorphic to the space of all linear
systems
$|D|$, where $D$ is a divisor with $c_1(\OO_X(2D-K))=-L$.
\end{pr}
\pf Use Theorem 4.2 and Bradlow's description of the moduli spaces of stable
pairs in the case of line bundles [B1].
\qed

We have already noticed (Remark 2.5) that in the case of a Riemannian
4-manifold
with nonnegative scalar curvature $s_g$, the Seiberg-Witten equations have only
reducible solutions. In the K\"ahler case, the same result can be obtained
under
the weaker assumption $\sigma_g \geq 0$ on the total scalar curvature.
\begin{co}
Let $(X,g)$ be a K\"ahler surface with nonnegative total scalar curvature
$\sigma_g$. Then all solutions of the Seiberg-Witten equations in rank 1 are
reducible. If moreover the surface has $K^2>0$, then for every almost
canonical class
$L$, the corresponding Seiberg-Witten equations are incompatible.
\end{co}
\pf The first assertion follows directly from the theorem, since the condition
$\sigma_g \geq 0$ is equivalent to $K\cup[\omega_g]\leq 0$. For the second
assertion, note that if $L$ is an almost canonical class, then $L^2=K^2>0$,
hence (regarded as line bundle) it cannot admit anti-selfdual connections.
\qed

\begin{re}
The Seiberg-Witten invariants associated to almost canonical classes are
well-defined for oriented, closed 4-manifolds $X$ satisfying $3\sigma+2e>0$.
\end{re}
\pf Recall that  if $L$ is an almost canonical class, then the expected
dimension of the moduli space of solutions of the perturbed Seiberg-Witten
equations [W, KM]corresponding to a $\Spin^c(4)$-structure of determinant $L$
is 0. Seiberg and Witten associate to every such class $L$ the number $n_L$
of points (counted with the correct signs [W]) of such a moduli space chosen to
be smooth and of the expected dimension. In the case $b_+\geq 2$, using
the same cobordism argument as in Donaldson theory, it follows that that
these numbers are well-defined, i.e. independent of the metric, provided the
moduli space has the expected dimension [KM]. The point is that the space of
$L$-good metrics [KM] (i.e. metrics with the property that the space of
harmonic anti-selfdual forms does not contain the harmonic representative of
$c_1^{\R}(L)$) is in this case path-connected. On the other hand, under the
assumption $3\sigma +2e>0$, it follows that
$L^2>0$ for any almost canonical class $L$, hence all metrics are $L$-good.
\qed

\begin{pr}
Let $(X,H_0)$ be a polarised surface with $K$ nef and big, and choose a
K\"ahler
metric $g$ with K\"ahler class $[\omega_g]=H_0+nK=:H$ for some $n\geq
KH_0$. Then
$\WW_X^g(L)$ is  empty for all almost canonical classes, except for $L=\pm K$,
when it consists of a  simple point.
\end{pr}
\pf Let $L$ be an almost canonical class  with $L H<0$.
Suppose $D$ is an effective divisor with $c_1(\OO_X(2D-K))=L$, so that
$D(D-K)=0$. Then $D^2=DK\geq 0$ since $K$ is nef. If $D^2$ were strictly
positive, the Hodge index theorem would give $(D-K)^2\leq 0$, i.e. $K^2\leq
D^2$.
But from $LH<0$ we get $0>(2D-K)(H_0+nK)=(2D-K)H_0+n(2D^2-K^2)\geq
(2D-K)H_0+n$,
which leads to the contradiction $n<(K-2D)H_0\leq KH_0$. Therefore
$D^2=DK=0$, so
that, again by the Hodge index theorem, $D$ must be numerically zero. Since $D$
is effective, it must be empty, and $L=-K$.

Replacing $L$ by $-L$ if $L$ is an almost
canonical class with $LH>0$, we find $L=K$ in this case. The corresponding
Seiberg-Witten moduli spaces are simple points in both cases, since
$H^2_{\bar\partial_A,\varphi}=H^1(\OO(D)|_D)=0$.
\qed

\begin{co}
There exists no orientation-preserving diffeomorphism between a rational
surface
and a minimal surface of general type.
\end{co}
\pf
Indeed, any rational surface $X$ admits a Hodge metric with
positive total scalar curvature [H2]. If $X$ was
orientation-preservingly diffeomorphic to a minimal surface of general
type, then
$K^2>0$, hence the Seiberg-Witten invariants are well defined (Remark 5.3),
and vanish by Corollary 5.2. Proposition 5.4 shows, however, that the
Seiberg-Witten invariants of a minimal surface of general type are non-trivial
for two almost canonical classes.
\qed
\vspace{3mm}\\

Witten has already proved [W] that for a
minimal surface of general type with $p_g>0$ ($b_+\geq 2$), the only almost
canonical classes which give non-trivial invariants are $K$ and $-K$. Their
proof
uses the moduli space of solutions of the perturbation of the Seiberg-Witten
equation with a holomorphic form. Proposition 5.4  shows that a stronger result
can be obtained with the  non-perturbed equations by choosing the  Hodge metric
$H=H_0+nK,\ n\gg 0$.

 For the proof of Corollary 5.5, we need in fact only the mod. 2 version of the
Seiberg-Witten invariants [KM2].

\newpage

\parindent0cm

\centerline {\Large {\bf Bibliography}}
\vspace{1cm}

[AHS] Atiyah M., Hitchin N. J., Singer I. M.: {\it Selfduality in
four-dimensional Riemannian geometry}, Proc. R. Lond. A. 362, 425-461 (1978)

[BPV] Barth, W., Peters, C., Van de Ven, A.: {\it Compact complex surfaces},
Springer Verlag (1984)

[B] Bertram, A.: {\it Stable pairs and stable parabolic pairs}, J. Alg.
Geometry
3, 703-724 (1994)

[B1] Bradlow, S. B.: {\it Vortices in holomorphic line bundles over closed
K\"ahler manifolds}, Comm. Math. Phys. 135, 1-17 (1990)

[B2] Bradlow, S. B.: {\it Special metrics and stability for holomorphic
bundles with global sections}, J. Diff. Geom. 33, 169-214 (1991)

[BD1] Bradlow, S. B.; Daskalopoulos, G.: {\it Moduli of stable pairs for
holomorphic bundles over Riemann surfaces I}, Intern. J. Math. 2, 477-513
(1991)

[BD2] Bradlow, S. B.; Daskalopoulos, G.: {\it Moduli of stable pairs for
holomorphic bundles over Riemann surfaces II}, Intern. J. Math. 4, 903-925
(1993)

[D] Donaldson, S.: {\it Anti-self-dual Yang-Mills connections over complex
algebraic surfaces and stable vector bundles}, Proc. London Math. Soc. 3,
1-26 (1985)

[DK] Donaldson, S.; Kronheimer, P.B.: {\it The Geometry of four-manifolds},
Oxford Science Publications (1990)

[FM] Friedman, R.; Morgan, J.W.: {\it Smooth 4-manifolds and Complex Surfaces},
Springer Verlag  3. Folge, Band 27 (1994)

[FQ]  Friedman, R.; Qin, Z.: {\it On complex surfaces diffeomorphic to
rational surfaces}, Preprint (1994)

[G1] Garcia-Prada, O.: {\it Dimensional reduction of stable bundles, vortices
and stable pairs}, Int. J. of Math. Vol. 5, No 1, 1-52 (1994)

[G2] Garcia-Prada, O.: {\it A direct existence proof for the vortex
equation over
a compact Riemann surface}, Bull. London Math. Soc., 26, 88-96 (1994)

[HH] Hirzebruch, F., Hopf H.: {\it Felder von Fl\"achenelementen in
4-dimensionalen 4-Mannigfaltigkeiten}, Math. Ann. 136 (1958)

[H1] Hitchin, N.: {\it  Harmonic spinors}, Adv. in Math. 14, 1-55 (1974)

[H2] Hitchin, N.: {\it  On the curvature of rational surfaces}, Proc. of Symp.
in Pure Math., Stanford, Vol. 27 (1975)

[HL1] Huybrechts, D.; Lehn, M.: {\it Stable pairs on curves and surfaces},
 J. Alg. Geometry  4, 67-104 (1995)

[HL2] Huybrechts, D.;  Lehn, M.: {\it Framed
modules and their moduli.}  Int. J. Math.
6,  297-324 (1995)

[JT] Jaffe, A., Taubes, C.: {\it Vortices and monopoles}, Boston, Birkh\"auser
(1980)

[K] Kobayashi, S.: {\it Differential geometry of complex vector bundles},
Princeton University Press (1987)

[KM1] Kronheimer, P.; Mrowka, T.: {\it Recurrence relations and asymptotics for
four-manifold invariants}, Bull. Amer. Math. Soc. 30, 215 (1994)

[KM2] Kronheimer, P.; Mrowka, T.: {\it The genus of embedded surfaces in the
projective plane}, Preprint (1994)

[OSS] Okonek, Ch.; Schneider, M.; Spindler, H: {\it Vector bundles  on complex
projective spaces}, Progress in Math. 3, Birkh\"auser, Boston (1980)

[Q] Qin, Z.: {\it Equivalence classes of polarizations and moduli spaces
of stable locally free rank-2 sheaves}, J. Diff. Geom. 37, No 2 397-416 (1994)

[S] Simpson, C. T.: {\it Constructing variations of Hodge structure using
Yang-Mills theory and applications to uniformization}, J. Amer. Math. Soc. 1
867-918 (1989)

[T]  Thaddeus, M.: {\it Stable pairs, linear systems and the Verlinde formula},
Invent. math. 117, 181-205 (1994)

[UY] Uhlenbeck,  K. K.; Yau,  S. T.: {\it On the existence of Hermitian
Yang-Mills connections in stable vector bundles}, Comm. Pure App. Math. 3,
257-293 (1986)

[W] Witten, E.: {\it Monopoles and four-manifolds}, Mathematical  Research
Letters 1,  769-796  (1994)
\vspace{2cm}\\
Authors addresses:\\
\\
Mathematisches Institut, Universit\"at Z\"urich,\\
Winterthurerstrasse 190, CH-8057 Z\"urich\\
e-mail:okonek@math.unizh.ch

\ \ \ \ \ \  \ \ \ teleman@math.unizh.ch

\end{document}